\pdfoutput=1

\RequirePackage{fix-cm}

\documentclass[leqno, fleqn, centertags, 12pt]{article}

\usepackage{latexsym}
\usepackage{amsmath}
\usepackage{amssymb}
\usepackage{verbatim}

\usepackage[T1]{fontenc}

\usepackage[upright, widespace]{fourier}

\usepackage{bm}

\usepackage{fullpage}

\usepackage{comment}

\newskip\nineskipamount \nineskipamount=9pt plus 0pt minus 0pt
\newskip\zeroskipamount \zeroskipamount=0pt plus 0pt minus 0pt
\usepackage[noorphans,vskip=\nineskipamount]{quoting}

\makeatletter
\renewcommand{\@makefntext}[1]{\vspace*{0.5ex}\parindent=0em
\hspace*{-0.4em}
\hbox to 0.4em{\hss\@makefnmark}\hspace*{0.4em}{#1}
}
\makeatother

\newcounter{mysectionnumber}
\newcommand{\mysection}[2]{\setcounter{footnote}{0}
\setcounter{myparnum}{0}
\refstepcounter{mysectionnumber}
\vspace{27pt}{\Large {\themysectionnumber.} {#1}}\label{#2}\vspace*{15pt}}

\newcommand{\myuppar}[1]{\vspace{\medskipamount}\textbf{#1}\hspace*{0.5em}}

\newcounter{myparnum}[mysectionnumber]
\renewcommand{\themyparnum}{\arabic{mysectionnumber}.\arabic{myparnum}}
\newcommand{\mypar}[2]{\refstepcounter{myparnum}{\vspace{\medskipamount}\textbf{{\themyparnum. #1}\label{#2}}\hspace{0.5em}}}

\newcounter{mylemmanum}[myparnum]
\newcounter{myappendnumber}
\newcounter{myaparnum}[myappendnumber]
\newcounter{myapparnum}[mysectionnumber]
\newcommand{\proof}{\vspace{\medskipamount}{\textbf{{\emph{Proof}.}}\hspace*{1em}}}

\newcommand{\eproof}{ $\blacksquare$}

\def\sss{\hspace{0.05em}\ }
\def\dss{\hspace{0.1em}\ }
\def\trs{\hspace{0.15em}\ }
\def\qss{\hspace{0.2em}\ }
\def\pss{\hspace{0.3em}\ }
\def\oss{\hspace{0.4em}\ }

\def\halfff{\hspace*{0.025em}}
\def\fff{\hspace*{0.05em}}
\def\dff{\hspace*{0.1em}}
\def\trf{\hspace*{0.15em}}
\def\qff{\hspace*{0.2em}}
\def\pff{\hspace*{0.3em}}
\def\off{\hspace*{0.4em}}

\newcommand{\nsp}{\hspace*{-0.1em}}
\newcommand{\nnsp}{\hspace*{-0.15em}}

\newcommand{\dnsp}{\hspace*{-0.2em}}

\renewcommand{\leq}{\leqslant}
\renewcommand{\geq}{\geqslant}

\newcommand{\zzz}{\mathbf{Z}}

\newcommand{\ccc}{\mathbf{C}}
\newcommand{\rrr}{\mathbf{R}}

\newcommand{\image}{\operatorname{Im}\trf}

\newcommand{\id}{\operatorname{id}}

\newcommand{\gr}{\operatorname{G{\fff}r}}

\newcommand{\norm}[1]{\|\qff #1 \qff\|}

\newcommand{\ttoo}{\hspace*{0.2em}\longrightarrow\hspace*{0.2em}}

\begin{document}

\setlength{\baselineskip}{12pt plus 0pt minus 0pt}
\setlength{\parskip}{12pt plus 0pt minus 0pt}
\setlength{\abovedisplayskip}{12pt plus 0pt minus 0pt}
\setlength{\belowdisplayskip}{12pt plus 0pt minus 0pt}

\newskip\smallskipamount \smallskipamount=3pt plus 0pt minus 0pt
\newskip\medskipamount   \medskipamount  =6pt plus 0pt minus 0pt
\newskip\bigskipamount   \bigskipamount =12pt plus 0pt minus 0pt

\author{Nikolai\qss V.\qss Ivanov}
\title{Spectral\qss sections\fff:\oss two\qss proofs\qss of\qss a\qss theorem\qss of\pss
Melrose--Piazza}
\date{}

\footnotetext{\hspace*{-0.65em}\copyright\oss 
Nikolai\qss V.\qss Ivanov,\oss 2021.\oss}

\footnotetext{\hspace*{-0.65em}The author\dss is\dss grateful\sss to\dss M.\dss Prokhorova\dss
for attracting\sss his attention\sss to\sss the work of\qss
Melrose\dss and\dss Piazza\qss \cite{mp1},\oss 
stimulating online discussions,\oss and writing\sss the paper\qss \cite{p3}.}

\maketitle

\renewcommand{\baselinestretch}{1}
\selectfont

\mysection{Introduction}{introduction}

\vspace{6pt}
Spectral\sss sections of\dss families of\dss
self-adjoint\trs Fredholm\dss operators were introduced\sss by\trs
Melrose\dss and\dss Piazza\qss \cite{mp1}\qss 
for\sss the needs of\dss index\sss theory.\oss
The basic result\sss about\sss spectral\sss sections\dss is\dss
a\sss theorem of\qss Melrose and\dss Piazza\trs to\sss the effect\sss that\sss
a family admits a spectral\sss section\dss if\trs and\dss only\trs if\trs
its analytic\sss index\sss vanishes.\oss
See\qss \cite{mp1},\qss Proposition\dss 1.\oss
Melrose\dss and\trs Piazza\qss \cite{mp1}\qss
gloss over\sss the definition of\dss the analytic\sss index and\sss 
the notion of\dss a\sss trivialization of\dss a\sss Hilbert\dss bundle\sss
implicitly\sss used\sss in\sss their\sss proof.\oss
Since passing\sss from one\sss trivializations of\dss 
a\dss Hilbert\dss bundle\sss to another rarely\sss preserves\sss
the norm continuity of\dss families of\trs Hilbert\sss space operators,\oss
the straightforward\sss interpretation of\trs Melrose--Piazza\dss proof\dss works only\sss
for\sss families of\dss operators in a fixed\dss Hilbert\dss space.\oss
The author\sss learned about\sss this\sss issue\sss from\dss M.\dss Prokhorova\qss \cite{p1}.

Recently\sss the author\sss proved\sss a\sss general\sss version of\dss this\dss
theorem of\trs Melrose--Piazza\dss as a byproduct\sss of\dss
a\sss theory developed\sss in\qss \cite{i1},\qss \cite{i2}.\oss
See\qss \cite{i2},\oss Corollary\qss 6.2.\oss
After\sss learning about\trs Melrose\dss approach\trs \cite{m}\trs 
to clarifying\qss \cite{mp1},\oss
the author\sss realized\sss that\sss a\sss less general\sss version of\qss Melrose--Piazza\dss theorem 
can be disentangled\sss from\sss the\sss theory of\pss \cite{i1},\qss \cite{i2}.\oss
This version\dss is\dss still\sss more general\sss than\sss the original\sss one\qss \cite{mp1}.\oss
An analysis of\dss the resulting\sss proof\dss led\sss to\dss
a fairly simple way\sss to prove\sss the original\trs Melrose--Piazza\dss theorem,\oss 
or,\oss rather,\oss its\qss ``axiomatic''\qss version.

The present\sss paper\dss is\dss devoted\sss to\sss 
the proofs\dss of\trs these\sss two versions
of\qss Melrose--Piazza\dss theorem.\oss
The ideas of\trs Atiyah--Singer\qss \cite{as}\qss play a\sss key\sss role,\oss
but,\oss in contrast\sss with\qss \cite{m}\qss and\qss \cite{as},\oss
compact\sss operators are not\sss used even\sss implicitly\sss
in\sss the proof\dss of\dss the first,\oss more general,\oss version.\oss
The proof\dss of\dss the second version\dss is\dss closer\sss 
to\sss the ideas of\trs Melrose\qss \cite{m}\qss and\sss uses 
compact\sss operators and\sss some ideas of\trs Atiyah--Segal\qss \cite{ase}.\oss
Both\sss proofs are presented as complements\sss to\qss \cite{mp1}.\oss
This forces us\sss to\sss deal only with compact\sss spaces of\dss parameters.\oss
If\dss one replaces\sss references\sss to\qss \cite{mp1}\qss
by\sss references\sss to\qss \cite{i2},\oss
the same proofs would\sss work for paracompact\sss spaces.\oss

Both\sss proofs depend on\sss a\sss theorem about\sss the existence of\dss
trivializations appropriately adapted\sss to families.\oss
See\dss Theorem\qss \ref{adapted}\qss below\sss for\sss the statement\sss 
and\qss \cite{i2},\oss Theorem\qss 4.6,\oss for a proof.\oss
The case of\dss triangulable bases\dss is\dss considered\sss in\qss \cite{i2},\oss
Theorem\qss 4.5,\oss and\dss is\dss independent\sss 
from\sss the rest\sss of\qss \cite{i2}.\oss
Also,\oss we refer\sss to\qss \cite{i2}\qss for\sss a detailed\sss proof\dss
of\trs Theorem\qss \ref{two-types-sections},\oss which\dss is\dss treated\sss
in\qss \cite{mp1}\qss as obvious.\oss
It\dss is\dss also independent\sss from\sss the rest\sss of\qss \cite{i2}.\oss

In\qss \cite{mp2}\qss Melrose\dss and\dss Piazza\dss
proved an odd $\zzz_{\dff 2}$\dnsp-graded version of\dss their\sss theorem.\oss
There\dss is\dss no doubt\sss that\sss the methods of\dss the present\sss paper\sss
work\sss in\sss this case also.\oss
Cf.\qss \cite{i1},\oss \cite{i2}.\oss

\mysection{Basic\qss definitions}{definitions}

\myuppar{The framework.}
Let\sss $H$\sss be a separable infinite dimensional\dss Hilbert\sss space over\sss $\ccc$\nnsp.\oss
Let\sss $\mathbb{H}$\sss be a\sss locally\sss trivial\dss Hilbert\sss bundle with\sss
a paracompact\dss base\sss $X$\sss and\sss fibers isomorphic\sss to\sss $H$\nnsp.\oss
We will\sss treat\sss $\mathbb{H}$\sss as a family\sss
$H_{\dff x}\dff,\qff x\qff \in\qff X$\sss 
of\trs Hilbert\sss spaces.\oss
We will\sss assume\sss that\sss the space\sss $X$\sss is\dss paracompact.\oss
Let\sss 
$A_{\dff x}\dff \colon\dff
H_{\dff x}\qff \ttoo\qff H_{\dff x}\dff,\pff 
x\qff \in\qff X$\sss 
be a family of\dss self-adjoint\trs Fredholm\dss operators.\oss
We will\sss assume\sss that\sss either\sss all\sss operators $A_{\dff x}$ are bounded,\oss
or\sss that\sss they are closed and densely defined.\oss
In\sss the first\sss case  
we will\sss assume\sss that\sss the family\sss
$A_{\dff x}\dff,\pff x\qff \in\qff X$\sss
is\dss continuous as a self-map of\dss the\sss total\sss space of $\mathbb{H}$\nnsp.\oss
In\sss the second case we will\sss assume\sss 
that\sss the\qss \emph{bounded\dss transform}\qss
$\gamma\trf(\trf A_{\dff x}\trf)\dff,\pff x\qff \in\qff X$\nnsp,\oss where\sss 
$\gamma\trf(\trf t\trf)
\off =\off 
t\trf(\trf 1\qff +\qff t^{\dff 2}\trf)^{\dff -\dff 1/2}$\dnsp,\oss
has\sss this property.\oss
By\sss a well\sss known\sss reason we assume\sss that\sss 
operators\sss $A_{\dff x}$\sss are neither 
essentially\sss positive,\oss nor\sss essentially negative.\oss

\myuppar{Enhanced operators.}
An\qss \emph{enhanced\qss (self-adjoint\trs Fredholm)\qss operator}\pss
is\dss a pair $(\trf A\fff,\qff \varepsilon\dff)$\dnsp,\pss
where\sss $A\dff \colon\dff H\qff \ttoo\qff H$\sss
is\dss a self-adjoint\trs Fredholm\dss operator,\pss
$\varepsilon\qff >\qff 0$\nnsp,\oss
such\sss that\sss
$\varepsilon\fff,\qff -\qff \varepsilon\qff \not\in\qff 
\sigma\trf(\trf A\trf)$\sss and\sss
the interval\sss 
$[\trf -\qff \varepsilon\fff,\qff \varepsilon\trf]$\sss
is\dss disjoint\sss from\sss the essential\sss spectrum of\sss $A$\nnsp.\oss
Then\sss the spectral\dss projection\sss
$P_{\dff [\dff -\dff \varepsilon\fff,\qff \varepsilon\dff]}\trf
(\trf A\trf)$\sss
has finitely dimensional\sss image.\oss
If\dss $(\trf A\fff,\qff \varepsilon\trf)$\sss is\dss an enhanced operator\sss
and\sss $A'\dff \colon\dff H\qff \ttoo\qff H$\sss is\dss a
self-adjoint\trs Fredholm\dss operator sufficiently close\sss to\sss $A$\sss
in\sss the norm\sss topology or\sss in\sss the uniform\sss resolvent\sss sense,\oss
 then\sss $(\trf A'\fff,\qff \varepsilon\trf)$\sss
is\dss also an enhanced operator.\oss
The spectral\sss projection\sss
$P_{\dff [\dff -\dff \varepsilon\fff,\qff \varepsilon\dff]}\trf
(\trf A'\trf)$ 
continuously depend on\sss $A'$ for\sss $A'$ 
sufficiently close\sss to\sss $A$\nnsp.\oss

\myuppar{Fredholm\dss families.}
The family\sss
$A_{\dff x}\fff,\qff x\qff \in\qff X$\sss
is\dss said\sss to be a\qss
\emph{Fredholm\dss family}\qss
if\dss all\sss operators $A_{\dff x}$ are\dss Fredholm\dss
and\sss for every\sss $z\qff \in\qff X$\sss
there exists\sss 
$\varepsilon\off =\off \varepsilon_{\dff z}\qff >\qff 0$\sss
and a neighborhood\sss $U_{\dff z}$\sss of\sss $z$
such\sss that\sss
$(\trf A_{\dff y}\dff,\qff \varepsilon\trf)$\sss
is\dss an enhanced\sss operator\sss 
for every $y\qff \in\qff U_{\dff z}$,\oss
the subspaces\vspace{1.5pt}\vspace{-0.125pt}
\[
\quad
V_{\fff y}
\off =\off
\image P_{\qff [\dff -\dff \varepsilon\fff,\qff \varepsilon\dff]}\trf(\trf A_{\dff y}\trf)
\off \subset\off 
H_{\dff y}
\]

\vspace{-12pt}\vspace{1.5pt}
continuously depend on\sss $y\qff \in\qff U_{\dff z}$,\oss
and\sss the operators\sss
$V_{\fff y}
\qff \ttoo\qff
V_{\fff y}$
induced\sss by\sss $A_{\dff y}$\sss
also continuously depend on\sss $y\qff \in\qff U_{\dff z}$.\oss
Two\sss families\sss
$A_{\dff x}\fff,\qff x\qff \in\qff X$\sss
and\dss
$A'_{\dff x}\fff,\qff x\qff \in\qff X$\sss
are said\sss to be\qss \emph{Fredholm\dss homotopic}\pss if\dss there exists\sss
a\dss homotopy\sss
$A_{\trf x\fff,\dff u}$,\qss 
$(\trf x\fff,\qff u\trf)
\qff \in\qff 
X\dff \times\dff [\dff 0\fff,\qff 1\dff]$\sss
which\dss is\dss Fredholm\dss as a family and\sss
such\sss that\sss
$A_{\trf x\fff,\dff 0}\off =\off A_{\dff x}$\sss
and\dss
$A_{\trf x\fff,\dff 1}\off =\off B_{\dff x}$\sss
for every\sss $x\qff \in\qff X$\nnsp.\oss

\myuppar{Strictly\trs Fredholm\dss families.}
Suppose\sss first\sss that\sss $\mathbb{H}$\sss is\dss the\sss trivial\sss bundle
with\sss the fiber\sss $H$\nnsp.\oss
We will\sss say\sss that\sss the family\sss
$A_{\dff x}\dff \colon\dff H\qff \ttoo\qff H$\nsp,\qss
$x\qff \in\qff X$\sss
is\qss \emph{strictly\trs Fredholm}\pss
if\dss it\dss is\trs Fredholm\dss and\sss
for every\sss $z\qff \in\qff X$\sss there exists\sss $\varepsilon\qff >\qff 0$\sss
and\dss a neighborhood\sss $U_{\fff z}$\sss of\sss $z$\sss
such\sss that\sss $(\trf A_{\dff y}\fff,\qff \varepsilon\trf)$\sss
is\dss an enhanced operator for every\sss 
$y\qff \in\qff U_{\fff z}$ and\sss the spectral\sss projection\sss
$P_{\dff [\trf \varepsilon\fff,\qff \infty\dff)}\trf
(\trf A_{\dff y}\trf)$\sss
continuously depends on\sss $y$\sss
in\sss the norm\sss topology\sss for\sss 
$y\qff \in\qff U_{\fff z}$\nsp.\oss

In\sss general,\pss 
$A_{\dff x}\dff \colon\dff H_{\dff x}\qff \ttoo\qff H_{\dff x}$\nsp,\qss
$x\qff \in\qff X$\sss
is\dss said\sss to be a\qss \emph{strictly\trs Fredholm\dss family}\pss
if\dss for every\sss $x\qff \in\qff X$\sss there exists\sss 
a neighborhood\sss $U_{\fff z}$\sss of\sss $z$\sss
and a\sss local\sss trivialization of\dss $\mathbb{H}$\sss
over\sss $U_{\fff z}$\sss turning\sss the restriction\sss
$A_{\dff x}\dff \colon\dff H_{\dff x}\qff \ttoo\qff H_{\dff x}$\nsp,\qss
$x\qff \in\qff U_{\fff z}$\sss
into a strictly\trs Fredholm\dss family\sss
in\sss the above sense.\oss
Such a\sss local\sss trivialization\dss is\dss said\sss to be\qss
\emph{strictly adapted}\pss to\sss the family\sss 
$A_{\dff x}\dff,\qff x\qff \in\qff X$\nnsp.\oss
A\sss trivialization of\sss $\mathbb{H}$\sss
is\dss said\sss to be\qss \emph{strictly adapted}\pss to\sss  
$A_{\dff x}\dff,\qff x\qff \in\qff X$\sss
if\dss for every\sss $z\qff \in\qff X$\sss there\sss exists\sss
a neighborhood\sss $U_{\fff z}$\sss of\sss $z$\sss such\sss that\sss
its restriction\sss to\sss $U_{\fff z}$\sss is\dss strictly adapted\sss to\sss
$A_{\dff x}\dff,\qff x\qff \in\qff X$\nnsp.

\myuppar{Fully\dss Fredholm\dss families.}
Suppose\sss that\sss 
$A_{\dff x}\dff \colon\dff H_{\dff x}\qff \ttoo\qff H_{\dff x}\dff,\qff
x\qff \in\qff X$\sss
is\dss a family of\dss bounded operators.\oss
We\sss say\sss that\sss it\dss is\qss \emph{fully\trs Fredholm}\pss
if\dss for every\sss $x\qff \in\qff X$\sss there exists\sss 
a neighborhood\sss $U_{\fff z}$\sss of\sss $z$\sss
and a\sss local\sss trivialization of\dss $\mathbb{H}$\sss
over\sss $U_{\fff z}$\sss turning\sss the restriction\sss
$A_{\dff x}\dff \colon\dff H_{\dff x}\qff \ttoo\qff H_{\dff x}$\nsp,\qss
$x\qff \in\qff U_{\fff z}$\sss
into a family continuous in\sss the norm\sss topology.\oss
Such a\sss local\sss trivialization\dss is\dss said\sss to be\qss
\emph{fully adapted}\pss to\sss  
$A_{\dff x}\dff,\qff x\qff \in\qff X$\nnsp.\oss
Clearly,\oss a\dss fully\trs Fredholm\dss family\dss is\dss strictly\trs Fredholm.\oss

If\dss 
$A_{\dff x}\dff \colon\dff H_{\dff x}\qff \ttoo\qff H_{\dff x}$\nsp,\qss
$x\qff \in\qff X$\sss
is\dss a family of\dss closed densely defined operators,\oss
we will\sss say\sss that\sss it\dss is\qss \emph{fully\trs Fredholm}\pss
if\dss the bounded\dss transform\sss
$\gamma\trf(\trf A_{\dff x}\trf)\dff,\pff x\qff \in\qff U_{\fff z}$\sss
is\dss fully\dss Fredholm.\oss
The\qss \emph{fully adapted}\pss ({\fff}local{\fff})\qss trivializations of\dss such families
are defined\sss in\sss the obvious way.\oss

\myuppar{Polarizations and\sss restricted\dss Grassmannians.}
A\qss \emph{polarization}\pss of\dss a\sss Hilbert\sss space\sss $K$\sss
is\dss a presentation of\sss $K$\sss as an orthogonal\sss direct\sss sum\dss
$K\off =\off K_{\dff -}\dff \oplus\dff K_{\dff +}$\dss
of\dss two closed\sss infinitely di\-men\-sion\-al\sss subspaces\sss
$K_{\dff -}\dff,\off K_{\dff +}$\nsp.\oss
A polarization\sss leads\sss to\sss the\qss
\emph{restricted\dss Grassmannian}\qss $\gr$\nnsp,\oss
the space of\dss subspaces\sss $L\qff \subset\qff K$\qss
\emph{commensurable}\pss with $K_{\dff -}$\nsp,\oss i.e.\qss
such\sss that\sss the intersection\dss 
$L\dff \cap\dff K_{\dff -}$\sss
is\dss closed and\sss has finite codimension\sss in\sss both\sss $L$\sss and\dss
$K_{\dff -}$\nsp.\oss
The\sss topology of\sss $\gr$\sss
is\dss defined\sss by\sss the norm\sss topology of\dss orthogonal\dss projections\sss
$K\qff \ttoo\qff L$\nnsp.\oss

\myuppar{Grassmannian\dss bundles and\sss weak\sss spectral\sss sections.}
Suppose\sss that\sss 
$A_{\dff x}\dff \colon\dff H_{\dff x}\trf \ttoo\qff H_{\dff x}$\nsp,\qss
$x\qff \in\qff X$\sss
is\dss a strictly\trs Fredholm\dss family.\oss
If\sss $x\qff \in\qff X$\sss and\sss
$(\trf A_{\dff x}\fff,\qff \varepsilon\trf)$\sss 
is\dss an enhanced operator,\oss then\vspace{3pt}
\[
\quad
H_{\dff x}
\off =\off
\image\dff
P_{\dff (\dff -\qff \infty\fff,\qff \varepsilon\trf]}\trf
(\trf A_{\dff x}\trf)
\qff \oplus\qss
\image\dff
P_{\dff [\trf \varepsilon\fff,\qff \infty\dff)}\trf
(\trf A_{\dff x}\trf)
\]

\vspace{-12pt}\vspace{3pt}
is\dss a polarization of\sss $H_{\dff x}$\nsp.\oss
This polarization\sss leads\sss to\sss a\sss
restricted\dss Grassmannian,\oss
which\sss we will\sss denote by\sss $\gr\trf(\trf x\trf)$\nnsp.\oss
Since\sss the family\sss $A_{\dff x}\dff,\qff x\qff \in\qff X$\sss
is\dss strictly\trs Fredholm,\oss
the family\sss of\dss restricted\dss Grassmannians\sss 
$\gr\trf(\trf x\trf)\dff,\qff
x\qff \in\qff X$\sss
forms a\sss locally\sss trivial\sss bundle\dss 
$\bm{\pi}\trf(\trf \mathbb{A}\trf)\dff \colon\dff
\gr\trf(\trf \mathbb{A}\trf)
\qff \ttoo\qff 
X$\nnsp.\oss
Continuous sections of\dss this bundle are called\qss
\emph{weak spectral\sss sections}\pss of\dss the family\sss
$A_{\dff x}\dff,\qff x\qff \in\qff X$\nnsp.

\myuppar{Discrete-spectrum\dss families and spectral\sss sections.}
The family\sss
$A_{\dff x}\dff \colon\dff H_{\dff x}\qff \ttoo\qff H_{\dff x}$\nsp,\qss
$x\qff \in\qff X$\sss
is\dss said\sss to be a\qss \emph{discrete-spectrum\dss family}\pss
if\dss for every\sss $\lambda\qff \in\qff \rrr$\sss the family\sss of\dss operators\sss 
$A_{\dff x}\qff -\pff \lambda$\nsp,\qss $x\qff \in\qff X$\sss
is\dss a\dss Fredholm\dss family.\oss
In\sss particular,\oss operators\sss $A_{\dff x}\qff -\pff \lambda$\sss
are\dss Fredholm\dss for every\sss $\lambda\qff \in\qff \rrr$\nnsp.\oss
The operators $A_{\dff x}$ in\dss such a family\sss have discrete spectrum
and cannot\sss be bounded.\oss
Since\sss they are self-adjoint,\oss they are necessarily closed and densely defined.\oss

Let\sss $A_{\dff x}\dff,\qff x\qff \in\qff X$\sss
be a strictly\trs Fredholm\dss and\sss discrete-spectrum\sss family.\oss
A\qss weak spectral\sss section\sss
$S\dff \colon\dff
X\qff \ttoo\qff \gr\trf(\trf \mathbb{A}\trf)$\sss 
is\dss said\sss to be a\qss\emph{spectral\sss section}\pss
if\trs \vspace{3pt}
\[
\quad
\image
P_{\qff [\trf r\dff(\dff x\trf)\fff,\pff \infty\qff)}\qff(\dff A_{\dff x}\dff)
\off \subset\off\dff
S\dff(\trf x\trf)
\off \subset\off\dff
\image
P_{\qff [\qff -\qff r\dff(\dff x\trf)\fff,\pff \infty\qff)}\qff(\dff A_{\dff x}\dff)
\pff
\]

\vspace{-12pt}\vspace{3pt}
for a continuous function
$r\dff \colon\dff X\qff \ttoo\qff \rrr_{\trf >\dff 0}$\sss 
and every\sss $x\qff \in\qff X$\nnsp.\oss
This\dss is\dss a generalization of\dss the notion of\dss spectral\sss sections 
introduced\sss by\trs Melrose\dss and\dss Piazza\qss \cite{mp1}.\oss

\myuppar{Discrete-spectrum\dss families and classical\sss operator\sss topologies.}
The material\sss of\dss this subsection\dss is\dss 
not\sss used\sss in\sss the rest\sss of\dss the paper.\oss

The discrete-spectrum\dss families were introduced\sss by\sss the author\qss \cite{i2}\qss
as a natural\sss analogue of\dss the notion of\dss a\dss Fredholm\dss family\sss
for families of\dss operators with discrete spectrum.\oss
At\sss the same\sss time\sss such families are exactly\sss the
families for which\sss the proof\dss of\qss Theorem\qss \ref{two-types-sections}\qss
below about\dss the existence of\dss spectral\sss sections works.\oss

Recently\dss M.\dss Prokhorova\qss \cite{p3}\qss related\sss the notion of\dss a discrete-spectrum\sss family 
with\sss classical\sss continuity\sss properties.\oss
Note\sss that\sss every operator in a discrete-spectrum\sss family\dss
is\dss a closed densely defined operator with compact\sss resolvent.\oss
Let\sss 
$A_{\dff x}\dff \colon\dff H\qff \ttoo\qff H$\nnsp,\qss
$x\qff \in\qff X$\sss
be a family\sss of\dss self-adjoint\dss operators
with compact\sss resolvent\sss in a fixed\dss Hilbert\sss space\sss $H$\nnsp.\oss
Then\sss $A_{\dff x}$\nsp,\qss $x\qff \in\qff X$\sss is\dss a discrete-spectrum\sss family\dss
if\trs and\dss only\trs if\sss
$A_{\dff x}$\nsp,\qss $x\qff \in\qff X$\sss
is\dss continuous in\sss the\sss topology of\dss convergence in\sss the norm\sss resolvent\sss sense.\oss
Also,\oss for discrete-spectrum\sss families every\sss strictly adapted\dss
local\dss trivialization\dss is\dss fully adapted.\oss
In\sss particular,\oss a discrete-spectrum\sss family\dss is\dss strictly\dss Fredholm\trs
if\trs and\dss only\trs if\dss it\dss is\dss fully\trs Fredholm.\oss
See\qss \cite{p3},\oss Theorems\qss 2\qss and\qss 3.

\mysection{Basic\qss theorems}{basic-theorems}

\myuppar{Classical\sss contractibility\sss theorems.}
Let\sss $K$\sss be a separable infinite dimensional\dss Hilbert\sss space.\oss
Let\sss us consider\sss the group of\dss isometries\sss
$K\qff \ttoo\qff K$\sss and denote by $U\dff(\trf K\trf)$ 
this group equipped with\sss the norm\sss topology.\oss
The group $U\dff(\trf K\trf)$\sss is\dss contractible by a\sss theorem of\trs
Kuiper\qss \cite{ku}.\oss

Another useful\sss topology\sss on\sss this group\dss is\dss the compact-open one.\oss
We will\sss denote by $\mathcal{U}\dff(\trf K\trf)$\sss this group equipped\sss with\sss
topology induced\sss from\sss the product\sss of\dss compact-open\sss topologies
by\sss the map\sss
$g
\off \longmapsto\off 
(\trf g\fff,\qff g^{\dff -\dff 1}\trf)$\nnsp.\oss 
Actually,\oss this\sss topology coincides with\sss the strong operator\sss topology,\oss
but\sss the author prefers\sss to ignore\sss this fact.\oss
The group $\mathcal{U}\dff(\trf K\trf)$\sss is\dss contractible by a\sss theorem of\trs
Atiyah\dss and\dss Segal\qss \cite{ase},\oss
who adapted an argument\sss of\qss 
Dixmier\dss and\dss Douady\qss \cite{dd}.\oss

Another\sss important\sss space\dss is\dss the space of\dss polarizations\sss
$H\off =\off H_{\dff -}\dff \oplus\dff H_{\dff +}$\sss
with\sss the\sss topology defined\sss by\sss the norm\sss topology of\dss
orthogonal\sss projections\sss $H\qff \ttoo\qff H_{\dff -}$\nsp.\oss
By an observation of\trs Atiyah\dss and\dss Singer\qss \cite{as},\oss
Kuiper's\trs theorem\sss implies\sss that\sss 
this space\dss is\dss also contractible.\oss

\mypar{Theorem.}{adapted}
\emph{If\trs
$A_{\dff x}\dff,\qff x\qff \in\qff X$\sss
is\dss a strictly\trs Fredholm\dss family,\oss
then\sss there exists a\sss trivialization of\trs $\mathbb{H}$\dss 
strictly\sss adapted\dss to\sss $A_{\dff x}\dff,\qff x\qff \in\qff X$\nnsp.\oss}

\proof
The proof\trs is\dss simpler\sss if\dss there exists a\sss triangulation of\dss
the space\sss $X$\nnsp,\oss perhaps infinite.\oss
In\sss this case one can\sss argue by\sss an\sss induction\sss by skeletons,\oss
using\sss at\sss each\sss step both\sss the contractibility of\dss
the space of\dss polarizations and\sss 
the contractibility of\dss the groups $\mathcal{U}\dff(\trf K\trf)$\nnsp.\oss
The\sss latter\dss is\dss applied\sss to\sss
$K\off =\off H_{\dff -}\dff,\pff H_{\dff +}$\sss
for polarizations\sss $H\off =\off H_{\dff -}\dff \oplus\dff H_{\dff +}$\nsp.\oss
The case of\dss triangulable space\sss $X$\sss is\dss sufficient\sss for applications.\oss
The general\sss case of\dss a paracompact\sss space\sss $X$\sss requires more sophisticated\sss
tools from\sss the homotopy\sss theory.\oss
See\qss \cite{i2},\oss Theorems\qss 4.5\qss and\qss 4.6.\oss  \eproof

\mypar{Theorem.}{two-types-sections}
\emph{If\trs
$A_{\dff x}\dff,\qff x\qff \in\qff X$\sss
is\dss a\sss discrete-spectrum\sss and\sss strictly\trs Fredholm\dss family,\oss
then every\sss weak\sss spectral\sss section\sss of\trs
$A_{\dff x}\dff,\qff x\qff \in\qff X$\sss 
is\dss homotopic\sss to
a spectral\sss section.\oss}

\proof
This\dss is\dss an explicit\sss form of\dss the\sss
last\sss paragraph\sss in\sss the proof\dss of\qss Proposition\trs 1\trs
of\pss Melrose\dss and\trs Piazza\qss \cite{mp1},\oss
who claim\sss that\sss a weak spectral\sss section can\sss be\sss
transformed\sss into a spectral\sss section\qss
\emph{``simply\dss by\sss smoothly\dss truncating\sss
the eigenfunction\sss expansion''.}\oss
See\qss \cite{i2},\oss Theorem\qss 6.1\qss for a detailed\sss geometric\sss version
of\dss this argument.\oss  \eproof

\mypar{Theorem.}{mp}
\emph{Let\sss
$A_{\dff x}\dff \colon\dff H\qff \ttoo\qff H$\nsp,\qss
$x\qff \in\qff X$\sss
be a norm continuous family of\qss Fredholm\dss self-adjoint\dss operators in a
fixed\dss Hilbert\sss space\sss $H$\nnsp.\oss
Then\sss 
$A_{\dff x}\dff,\qff x\qff \in\qff X$\sss 
admits a weak\sss spectral\sss section\dss if\trs and\dss only\trs if\trs it\dss
is\dss homotopic\sss in\sss the class of\dss such\sss families\sss
to a family\sss of\dss invertible operators.\oss}

\proof
For compact\sss $X$\sss the proof\dss is\dss contained\sss in\qss \cite{mp1}.\oss 
See\qss \cite{mp1},\oss the proof\dss of\qss Proposition\qss 1.\oss
For\sss general\sss paracompact\sss $X$\sss this follows from\qss \cite{i2}.\oss 
We omit\sss the details.\oss \eproof

\myuppar{Remark.}
Suppose\sss that\sss 
$A_{\dff x}\fff \colon\fff H\dff \ttoo\dff H$\nnsp,\dss
$x\qff \in\qff X$\sss
is\dss a\sss family of\qss (closed densely defined{\halfff})\qss self-adjoint\sss operators.\oss
The results of\trs Prokhorova\qss \cite{p3}\qss mentioned at\sss the end of\trs
Section\qss \ref{definitions}\qss imply\sss that\sss
$A_{\dff x}$\nsp,\dss $x\qff \in\qff X$\sss
is\dss a discrete-spectrum and strictly\dss Fredholm\dss family\dss
if\dss and only\dss if\dss $A_{\dff x}$\sss are
operators with compact\sss resolvent\sss
and\sss the family of\dss bounded\sss transforms\sss $\gamma\trf(\trf A_{\dff x}\trf)$\nsp,\dss
$x\qff \in\qff X$\sss is\dss norm\sss continuous.\oss
In view of\dss this equivalence,\oss 
the analogues of\qss Theorems\qss \ref{two-types-sections}\qss
and\qss \ref{mp}\qss for such families with paracompact\sss $X$\sss follow\sss
from\sss the results of\trs Prokhorova\qss \cite{p2}.\oss
See\qss \cite{p2},\oss Theorem\qss 4.4\qss (note\sss that\sss every weak spectral\sss section\dss
is\dss a generalized spectral\sss section).\oss

\mysection{The\qss first\qss proof}{first-proof}

\myuppar{Finite-polarized\sss replacements.}
We will\sss call\sss a self-adjoint\sss operator\sss
$A\dff \colon\dff K\qff \ttoo\qff K$\sss in a\sss Hilbert\sss space\sss $K$\dss
\emph{finite-polarized}\oss if\dss
$\norm{A}\off =\off 1$\nnsp,\oss
the essential\sss spectrum of\sss $A$\sss consists of\dss two points\sss
$-\qff 1\fff,\qff 1$\nnsp,\oss
and\dss the spectral\sss projection\sss
$P_{\trf (\dff -\qff 1\fff,\qff 1\trf)}\trf(\trf A\trf)$\sss
is\dss an operator of\dss finite rank.\oss
If\dss we omit\sss the\sss last\sss property,\oss 
we will\sss get\sss exactly\sss the operators from\sss the space\sss $\hat{F}_{\dff *}$
from\qss \cite{as},\oss Section\qss 2.\oss

Suppose\sss that\sss the family\sss 
$A_{\dff x}\dff,\qff x\qff \in\qff X$\sss
is\dss Fredholm.\oss
We would\sss like\sss to replace\sss it\sss
by a family\sss of\dss finite-polarized operators\sss
$A'_{\dff x}\dff \colon\dff 
H_{\dff x}\qff \ttoo\qff H_{\dff x}\dff,\qff 
x\qff \in\qff X$\sss
without\sss affecting\sss analytic\sss index\sss
and\sss weak\sss spectral\sss sections.\oss
This can\sss be done by\sss a spectral\sss deformation\sss similar\sss to
one used\sss in\qss \cite{as}.\oss 
Let\sss us choose\sss for each\sss $x\qff \in\qff X$\sss
a neighborhood\sss $U_{\dff x}$\sss of\sss $x$
and a number\sss $\varepsilon_{\dff x}\qff \in\qff (\dff 0\fff,\qff 1\dff)$\sss
such\sss that\sss the properties from\sss the definition of\trs
Fredholm\dss families hold.\oss
Since\sss $X$\sss is\dss paracompact,\oss
we can assume\sss that\sss for some\sss $\Sigma\qff \subset\qff X$\sss
the family\sss $U_{\fff a}\dff,\qff a\qff \in\qff \Sigma$\sss
is\dss a\sss locally\sss finite covering of\sss $X$\sss and\sss
there exists a partition of\dss unity\sss subordinated\sss to\sss this covering.\oss
Let\sss 
$r_{\fff a}\dff \colon\dff X\qff \ttoo\qff \rrr_{\qff \geq\dff 0}$\sss 
be\sss the function from\dss this partition of\dss unity 
corresponding\sss to\sss $a\qff \in\qff \Sigma$\nnsp.\oss
Let\vspace{1.5pt}
\[
\quad
r\trf(\trf x\trf)
\off =\off
\sum\nolimits_{\qff a\qff \in\qff \Sigma}\qff r_{\fff a}\dff(\trf x\trf)\qff \varepsilon_{\dff a}
\pff.
\]

\vspace{-12pt}\vspace{1.5pt}
Then $r\trf(\trf x\trf)\qff \in\qff (\dff 0\fff,\qff 1\dff)$ and\sss
$r\trf(\trf x\trf)\qff \leq\qff \max\qff \varepsilon_{\dff a}$\nsp,\oss
where\sss the maximum\dss is\dss over\sss $a$\sss 
such\sss that\sss $x\qff \in\qff U_{\fff a}$\nsp,\oss
for every\sss $x\qff \in\qff X$\nnsp.\oss
It\sss follows\sss that\sss for every\sss $x\qff \in\qff X$\sss
the essential\sss spectrum of\sss $A_{\dff x}$\sss
is\dss disjoint\sss from\sss
$(\trf -\qff r\trf(\trf x\trf)\fff,\qff r\trf(\trf x\trf)\qff)$\nnsp.\oss
For\sss $r\qff >\qff 0$\sss let\sss 
$\chi_{\dff r}\dff \colon\dff \rrr\qff \ttoo\qff \rrr$\sss
be\sss an odd\sss increasing\sss function\sss such\sss that\vspace{1.5pt}
\[
\quad
\chi_{\dff r}\dff(\trf u\trf)
\off =\off
u
\quad
\mbox{for}\quad\qff
0\qff \leq\qff u\qff \leq\qff r/2\qff,
\quad
\mbox{and}\quad
\chi_{\dff r}\dff(\trf u\trf)
\off =\off
1
\quad
\mbox{for}\quad
u\qff \geq\qff r\qff.
\]

\vspace{-12pt}\vspace{1.5pt}
We assume\sss that\sss $\chi_{\dff r}$ continuously depends on $r$ in\sss
the $\sup$\dnsp-norm\sss topology.\oss
For\sss $x\qff \in\qff X$\sss let\vspace{1.5pt}
\[
\quad
A'_{\dff x}
\off =\off
\chi_{\dff r\trf(\trf x\trf)}\trf\bigl(\trf A_{\dff x}\trf\bigr)
\qff \colon\qff 
H_{\dff x}\qff \ttoo\qff H_{\dff x}
\pff
\]

\vspace{-12pt}\vspace{1.5pt}
The operators\sss $A'_{\dff x}$\sss are finite-polarized and\sss
$A'_{\dff x}\dff,\pff x\qff \in\qff X$\sss is\dss a\dss Fredholm\dss family.\oss
The family\sss $A'_{\dff x}\dff,\qff x\qff \in\qff X$\sss
is\dss our\qss \emph{finite-polarized\sss replacement}\pss
of\sss $A_{\dff x}\dff,\qff x\qff \in\qff X$\sss.\oss

\mypar{Lemma.}{strictly-f-replacement}
\emph{If\dss the family\dss 
$A_{\dff x}\dff,\qff x\qff \in\qff X$\sss is\dss strictly\dss Fredholm,\oss
then\sss $A'_{\dff x}\dff,\qff x\qff \in\qff X$\sss is\dss also strictly\dss Fredholm,\oss
and\sss these\sss two families have\sss the same weak spectral\sss sections.\oss}

\proof
Clearly,\oss if\sss
$0\qff <\qff \varepsilon\qff <\qff r\trf(\trf x\trf)/2$\nnsp,\oss
then\vspace{0.75pt}
\[
\quad
P_{\qff [\dff \varepsilon\fff,\qff \infty\dff)}\qff
\left(\trf A_{\dff x}\trf\right)
\off =\off\dss
P_{\qff [\dff \varepsilon\fff,\qff \infty\dff)}\qff
\left(\trf A'_{\dff x}\trf\right)
\pff.
\]

\vspace{-12pt}\vspace{0.75pt}
It\sss follows\sss that\dss if\dss 
$A_{\dff x}\dff,\qff x\qff \in\qff X$\sss is\dss strictly\dss Fredholm,\oss
then\sss $A'_{\dff x}\dff,\qff x\qff \in\qff X$\sss is\dss also strictly\dss Fredholm.\oss
This also implies\sss that\sss the bundles\sss
$\bm{\pi}\trf(\trf \mathbb{A}\trf)\dff \colon\dff
\gr\trf(\trf \mathbb{A}\trf)
\qff \ttoo\qff 
X$\sss
and\sss
$\bm{\pi}\trf(\trf \mathbb{A}'\trf)\dff \colon\dff
\gr\trf(\trf \mathbb{A}'\trf)
\qff \ttoo\qff 
X$\sss
corresponding\sss to\sss these\sss two families
are equal\qss (not\sss only\sss isomorphic,\oss but\sss equal\fff).\oss
Therefore\sss these bundles have\sss the same sections,\oss
i.e.\qss weak\sss spectral\sss sections are\sss the same.\oss  \eproof

\mypar{Lemma.}{replacements-fully-F}
\emph{Every\sss local\dss trivialization\sss strictly adapted\sss to\sss the family\sss
$A_{\dff x}\dff,\qff x\qff \in\qff X$\sss 
is\dss fully\sss adapted\dss to\sss the family\sss 
$A'_{\dff x}\dff,\qff x\qff \in\qff X$\nnsp.\oss 
In\sss particular,\oss
if\qss the family\dss 
$A_{\dff x}\dff,\qff x\qff \in\qff X$\sss is\dss strictly\dss Fredholm,\oss
then\sss the family\sss 
$A'_{\dff x}\dff,\qff x\qff \in\qff X$\sss
is\qss fully\dss Fredholm.\oss}

\proof
Let\sss $z\qff \in\qff X$\nnsp.\oss 
There exist\sss $\varepsilon\qff >\qff 0$\sss
and a neighborhood\sss $U_{\dff z}$\sss of\sss $z$\sss such\sss that\sss
$(\trf A_{\dff y}\dff,\qff \varepsilon\trf)$\sss
is\dss an enhanced operator\sss and\sss
$\varepsilon\qff <\qff r\trf(\trf y\trf)/2$\dss
for every\sss $y\qff \in\qff U_{\dff z}$\nsp.\oss
If\sss $U$\sss is\dss sufficiently small,\oss then
a strictly adapted\sss trivialization over\sss $U$\sss turns\sss the family\vspace{1.5pt}
\[
\quad
P_{\qff [\dff \varepsilon\fff,\qff \infty\dff)}\qff
\left(\trf A_{\dff y}\trf\right)
\off -\off\dff
P_{\qff (\dff -\qff \infty\fff,\qff \varepsilon\trf]}\qff
\left(\trf A_{\dff y}\trf\right)\dff,\quad
y\qff \in\qff U
\]

\vspace{-12pt}\vspace{1.5pt}
into a norm-continuous\sss one.\oss
But\sss $A'_{\dff y}\dff,\qff u\qff \in\qff U$\sss
differs from\sss it\sss by\sss a norm continuous family of\dss operators of\dss finite rank.\oss
This implies\sss that\sss the same\sss local\dss trivialization\sss
turns\sss the family\sss $A'_{\dff y}\dff,\qff u\qff \in\qff U$\sss
into a norm continuous one.\oss  \eproof

\mypar{Lemma.}{fully-F-homotopy}
\emph{If\dss the family\dss 
$A_{\dff x}\dff,\qff x\qff \in\qff X$\sss is\dss fully\dss Fredholm,\oss
then\sss there exists a\sss fully\trs Fredholm\dss homotopy\sss
between\dss $A_{\dff x}\dff,\qff x\qff \in\qff X$\dss
and\trs $A'_{\dff x}\dff,\qff x\qff \in\qff X$\nnsp.\oss}

\proof
It\sss is\dss sufficient\sss to consider\sss the case when\sss
$A_{\dff x}\dff,\qff x\qff \in\qff X$\sss is\dss a family of\dss bounded operators,\oss
since\sss the case of\dss closed densely defined operators reduces\sss to\sss it\sss
by applying\sss the bounded\sss transform\sss to both\sss families.\oss
The\sss linear homotopies between\sss the identity\sss 
$\id\dff \colon\dff \rrr\qff \ttoo\qff \rrr$\sss
and\sss the functions\sss
$\chi_{\dff r\trf(\trf x\trf)}$\sss
define a homotopy\sss between\sss
$A_{\dff x}\dff,\qff x\qff \in\qff X$\dss
and\trs $A'_{\dff x}\dff,\qff x\qff \in\qff X$\nnsp.\oss
If\sss $A_{\dff x}\dff,\qff x\qff \in\qff X$\sss 
is\dss fully\dss Fredholm,\oss
then\sss this homotopy\dss is\dss also fully\trs Fredholm.\oss  \eproof

\myuppar{Atiyah--Singer\dss approach\sss to\sss the analytic\sss index.}
Let\sss
$A_{\dff x}\dff \colon\dff H\qff \ttoo\qff H$\nsp,\qss
$x\qff \in\qff X$\sss
be a\dss Fredholm\dss family of\dss bounded operators in a fixed\dss Hilbert\sss space $H$\nnsp.\oss
Let\sss us consider\sss the family\sss of\trs Fredholm\dss operators\qss
(not\sss assumed\sss to be self-adjoint\fff)\qss 
defined\sss by\sss the formula\vspace{1.5pt}
\begin{equation}
\label{as-formula}
\quad
B_{\dff x\fff,\qff t}
\off =\off
\id_{\trf H}\qff \cos t
\pff +\pff
i\trf A_{\dff x}\qff \sin t
\qff,
\quad\off
x\qff \in\qff X\dff,\off\qff
t\qff \in\qff [\trf 0\fff,\qff \pi\trf]
\pff,
\end{equation}

\vspace{-12pt}\vspace{1.5pt}
where\sss $\id_{\trf H}$\sss is\dss the identity operator\sss
$H\qff \ttoo\qff H$\sss and\dss
$i\off =\off \sqrt{\dff -\qff 1}$\nnsp.\oss
The idea\dss is\dss to define\sss the analytic\sss index of\dss
$A_{\dff x}\dff,\qff x\qff \in\qff X$\sss
as\sss the analytic\sss index of\dss the family\sss
$B_{\dff x\fff,\qff t}\dff,\off 
(\trf x\fff,\qff t\trf)
\qff \in\qff
X\dff \times\dss [\trf 0\fff,\qff \pi\trf]$\nnsp.\oss
The\sss latter\sss should\sss be considered\sss relatively\sss
to\sss $X\dff \times\dff 0\qff \cup\qff X\dff \times\dff \pi$\nnsp,\oss
reflecting\sss the fact\sss that\sss
$B_{\dff x\fff,\qff 0}\off =\off \id_{\trf H}$\sss
and\sss $B_{\dff x\fff,\qff \pi}\off =\off -\qff \id_{\trf H}$\sss
for every\sss $x\qff \in\qff X$\nnsp.\oss
More precisely,\oss the analytic\sss index\dss is\dss either an element\sss
of\sss $K\trf(\trf \Sigma\dff X\trf)$\nnsp,\oss where\sss $\Sigma\dff X$\sss
is\dss the suspension of\sss $X$\sss 
or\sss a homotopy class of\dss a map\sss from\sss
$X\dff \times\dss [\trf 0\fff,\qff \pi\trf]$\sss
to an appropriate classifying space with respect\sss to homotopies
fixed on\sss $X\dff \times\dff 0$\sss and\sss $X\dff \times\dff \pi$\nnsp.\oss

If\sss $A_{\dff x}\dff,\qff x\qff \in\qff X$\sss is\dss fully\dss Fredholm,\oss
then\sss the family\qss (\ref{as-formula})\qss is\qss
\emph{fully\trs Fredholm}\pss in\sss the sense\sss that\sss
over an open neighborhood of\dss every\sss
$(\trf x\fff,\qff t\trf)
\qff \in\qff
X\dff \times\dff [\trf 0\fff,\qff \pi\trf]$\sss
there\sss exists a\sss trivialization of\sss 
$\mathbb{H}\dff \times\dff [\trf 0\fff,\qff \pi\trf]$\sss
turning\sss this family\sss into a norm-continuous family.\oss
The definition of\dss the analytic\sss index\sss
in\sss the non-self-adjoint\sss case\qss \cite{a},\oss \cite{as4}\qss
trivially\sss extends\sss to such\sss families,\oss
at\sss least\sss for\sss compact\sss $X$\nnsp.\oss
Hence we can extend\sss the definition of\dss the analytic\sss index\sss
to fully\dss Fredholm\dss families of\dss self-adjoint\sss operators.\oss
Clearly,\oss it\dss is\dss invariant\sss under\dss fully\trs Fredholm\dss homotopies.\oss
If\sss $X$\sss is\dss only\sss paracompact,\oss one needs\sss to use\dss
Segal's\dss definition\qss \cite{sfc}\qss of\dss the analytic\sss index 
and more\sss work\dss is\dss required.\oss
See\qss \cite{i2},\oss Sections\qss 7\qss and\qss 8,\oss
for\sss this case.\oss

It\dss is\dss easy\sss to see\sss that\sss the bounded\sss transform\sss
$\gamma\trf(\trf A_{\dff x}\trf)\dff,\pff x\qff \in\qff X$\sss
has\sss the same analytic\sss index as\sss
$A_{\dff x}\dff,\qff x\qff \in\qff X$\nnsp.\oss
This allows\sss to extend\sss the above definition of\dss the analytic index\sss
to families of\dss closed densely defined self-adjoint\trs Fredholm\dss operators in\sss $H$\nnsp.\oss
Namely,\oss given such a family\sss $A_{\dff x}\dff,\qff x\qff \in\qff X$\nnsp,\oss
one defines its index as\sss the index of\dss the bounded\dss transform\sss
$\gamma\trf(\trf A_{\dff x}\trf)\dff,\pff x\qff \in\qff X$\nnsp.\oss

Theorem\qss \ref{adapted}\qss allows\sss to extend\sss this definition\sss
to strictly\trs Fredholm\dss families\dss
of\dss operators in\sss the fibers of\dss a\dss Hilbert\dss bundle.\oss
The index does not\sss depend on\sss the choice of\dss trivialization\sss
by\sss the relative version of\qss Theorem\qss \ref{adapted}.\oss
See\qss \cite{i2},\oss Theorem\qss 4.7\qss for\sss the\sss latter.\oss

\mypar{Theorem}{index-wss}
\emph{Suppose\sss that\sss the family\dss 
$A_{\dff x}\dff,\qff x\qff \in\qff X$\sss is\dss fully\dss Fredholm.\oss
A weak\sss spectral\sss section\sss for\sss the family\sss
$A_{\dff x}\dff,\qff x\qff \in\qff X$\sss exists\dss if\trs and\dss only\trs if\trs
the analytic\sss index of\qss
$A_{\dff x}\dff,\qff x\qff \in\qff X$\sss vanishes.\oss}

\proof
We will\sss limit\sss ourselves by\sss the case of\dss compact\sss $X$\nnsp.\oss
Suppose\sss that\sss there exists a weak\sss spectral\sss section.\oss
Then\sss the arguments\sss in\sss the first\sss part\sss of\dss the proof\dss
of\qss Proposition\qss 1\qss in\qss \cite{mp1}\qss
together\sss with\sss the principle of\dss uniform\sss boundedness show\sss that\sss
there exists a\sss fully\trs Fredholm\dss homotopy\sss between\sss
$A_{\dff x}\dff,\qff x\qff \in\qff X$\sss and\sss a family\sss of\dss
invertible operators.\oss
Hence we may assume\sss that\sss operators\sss 
$A_{\dff x}\dff,\qff x\qff \in\qff X$\sss
are invertible.\oss
Then\sss all\sss operators\qss (\ref{as-formula})\qss are also invertible.\oss 
In\sss the non-self-adjoint\sss case\sss the analytic\sss index of\dss
families of\dss invertible operators vanishes essentially\sss by\sss the definition.\oss
The\qss ``only\trs if''\qss part\sss follows.\oss

Suppose now\sss that\sss the analytic\sss index of\sss
$A_{\dff x}\dff,\qff x\qff \in\qff X$\sss vanishes.\oss
Lemma\qss \ref{fully-F-homotopy}\qss implies\sss that\sss
the analytic\sss index of\dss our\sss finite-polarized\sss replacement\sss
$A'_{\dff x}\dff,\qff x\qff \in\qff X$\sss also vanishes.\oss
By\trs Theorem\qss \ref{adapted}\qss there exists a strictly adapted\sss to\sss
$A'_{\dff x}\dff,\qff x\qff \in\qff X$\sss trivialization of\sss $\mathbb{H}$\nnsp.\oss
By\trs Lemma\qss \ref{replacements-fully-F}\qss such\sss trivialization\sss
turns\sss the family\sss $A'_{\dff x}\dff,\qff x\qff \in\qff X$\sss
into a norm continuous one.\oss
The analytic index of\dss a norm continuous family vanishes\sss
if\trs and\dss only\trs if\dss it\dss is\dss homotopic\sss to a constant\sss family.\oss
Since a self-adjoint\sss operator can be deformed\sss to an\sss invertible
self-adjoint\sss operator,\oss in\sss this case\sss
$A'_{\dff x}\dff,\qff x\qff \in\qff X$\sss
is\dss homotopic\sss to a family of\dss invertible operators.\oss
By\trs Theorem\qss \ref{mp}\qss this implies\sss that\sss
there exists a\sss weak\sss spectral\sss section of\sss
$A'_{\dff x}\dff,\qff x\qff \in\qff X$\nnsp.\oss
Lemma\qss \ref{strictly-f-replacement}\qss implies\sss that\sss
this weak\sss spectral\sss section\dss is\dss also a weak\sss
spectral\sss section\sss for\sss 
$A_{\dff x}\dff,\qff x\qff \in\qff X$\nnsp.\oss
This proves\sss the\qss ``if''\qss part.\oss  \eproof

\mypar{Corollary.}{index-ss}
\emph{Suppose\sss that\dss 
$A_{\dff x}\dff,\qff x\qff \in\qff X$\sss is\dss 
a discrete-spectrum\sss and\dss fully\dss Fredholm\dss family.\oss
A spectral\sss section\sss for\sss 
$A_{\dff x}\dff,\qff x\qff \in\qff X$\sss 
exists\dss if\trs and\dss only\trs if\trs the analytic\sss index of\qss
$A_{\dff x}\dff,\qff x\qff \in\qff X$\sss vanishes.\oss}

\proof
It\dss is\dss sufficient\sss to combine\trs
Theorems\qss \ref{index-wss}\qss and\qss \ref{two-types-sections}.\oss  \eproof

\myuppar{Remark.}
In\trs Theorem\qss \ref{index-wss}\qss and\dss Corollary\qss \ref{index-ss}\qss
it\dss is\dss sufficient\sss to assume\sss that\sss
$A_{\dff x}\dff,\qff x\qff \in\qff X$\sss is\dss a strictly\trs Fredholm\dss family.\oss
This follows\sss from\qss \cite{i2},\oss Theorems\qss 5.2\qss and\dss
Corollary\qss 6.2,\oss if\dss one\sss takes into account\sss the equivalence of\dss
the definition of\dss the analytic index\sss suggested\sss in\qss \cite{i2}\qss
and\sss the classical\sss definition.\oss
This equivalence\dss is\dss proved\sss in\qss \cite{i2},\oss Theorem\qss 8.5.\oss
The proofs automatically\sss work\sss for paracompact\sss $X$\nnsp.\oss
The proof\dss of\dss the equivalence of\dss two definitions
of\dss the analytic\sss index\dss is\dss simpler\sss for compact\sss $X$\nnsp,\oss
where\sss in\sss the non-self-adjoint\sss case 
one can use\dss Atiyah's\dss definition\qss \cite{a}.\oss 
For\sss paracompact\sss $X$\sss one needs\sss to use\dss
Segal's\dss definition\qss \cite{sfc}.\oss
The resulting proof\dss of\dss
the\qss ``only\trs if''\qss parts depends on\sss the\sss theory developed\sss in\qss \cite{i1}.

\myuppar{Families of\dss elliptic operators.}
Let $\mathbb{M}$ be a\sss locally\sss trivial\sss bundle over $X$
with closed\sss manifolds as fibers.\oss
Let\sss us\sss consider
a continuous\sss family\sss of\dss elliptic pseudo-differential\sss operators
of\dss order\sss $0$\sss
acting on\sss fibers of\sss $\mathbb{M}$\nnsp.\oss  
It\sss defines a family of\dss bounded operators\sss acting in\sss 
fibers of\dss a\dss Hilbert\dss bundle\sss $\mathbb{H}$\sss over\sss $X$\nnsp.\oss
Classical\sss results of\trs Seeley\qss \cite{see}\qss and\dss Atiyah\dss and\dss Singer\qss \cite{as4}\qss
show\sss that\sss the\sss latter\sss family\dss is\dss fully\trs Fredholm.\oss
Hence\trs Theorem\qss \ref{index-wss}\qss 
applies\sss to such\sss families.\oss

Suppose now\sss that\sss we are given a
family of\dss elliptic self-adjoint\sss differential\sss operators of\dss order $1$
acting on\sss fibers of\sss $\mathbb{M}$\nnsp.\oss 
This\dss is\dss the class of\dss families considered\sss 
in\qss \cite{mp1},\oss Proposition\qss 1.\oss
Such a family defines a family\sss 
$A_{\dff x}\dff,\qff x\qff \in\qff X$\sss 
of\dss closed densely defined operators\sss acting in\sss 
fibers of\dss a\dss Hilbert\dss bundle\sss $\mathbb{H}$\sss over\sss $X$\nnsp.\oss
It\dss is\dss well\sss known\sss that\sss
$A_{\dff x}\dff,\qff x\qff \in\qff X$\sss
has\sss the properties defining\sss
discrete-spectrum\dss fully\trs Fredholm\dss families.\oss
Hence\dss Corollary\qss \ref{index-ss}\qss
applies\sss to such\sss families.\oss

\mysection{The second\dss proof}{second-proof}

\myuppar{Compactly-polarized\sss operators.}
Let\sss us\sss call\sss a self-adjoint\sss operator\sss
$A\dff \colon\dff K\qff \ttoo\qff K$\sss in a\sss Hilbert\sss space\sss $K$\dss
\emph{compactly-polarized}\oss if\dss
$\norm{A}\off =\off 1$\sss and\dss
the essential\sss spectrum of\sss $A$\sss consists of\dss two points\sss
$-\qff 1\fff,\qff 1$\nnsp,\oss
i.e.\qss $A$\sss belongs\sss to\dss Atiyah--Singer\qss \cite{as}\qss 
space\sss $\hat{F}_{*}$\nsp.\oss
Such operator\sss $A$\sss is\dss automatically\trs Fredholm.\oss
One can also define compactly-polarized operators as essentially unitary operators with\sss the norm $1$\nnsp.\oss
The\sss term\qss \emph{compactly-polarized}\oss is\dss intended\dss to stress\sss the obvious
analogy\sss with\sss finitely-polarized\sss operators from\dss Section\qss \ref{first-proof}.\oss

We claim\sss that\sss every compactly-polarized operator $A$\sss
has\sss the form\sss $Q\qff +\qff k$\nnsp,\oss
where $Q$\sss is\dss a unitary self-adjoint\sss operator and $k$\sss
is\dss a compact\sss self-adjoint\sss operator.\oss
Of\dss course,\oss this\dss is\dss well\sss known,\oss
but\sss we need some notations from\sss the proof.\oss
Let\sss $A$\sss be a compactly-polarized operator and\sss let\sss
$\varepsilon\qff \in\qff (\dff 0\fff,\qff 1\dff)$\sss be such\sss that\sss
$(\trf A\fff,\qff \varepsilon\trf)$\sss
is\dss an enhanced operator.\oss
Let\vspace{3pt}
\[
\quad
Q
\off =\off
Q_{\qff \varepsilon}\dff(\trf A\trf)
\off =\off 
P_{\qff [\dff \varepsilon\fff,\qff \infty\dff)}\trf(\trf A\trf)
\pff -\pff
P_{\qff (\dff -\qff \infty\fff,\qff \varepsilon\trf]}\trf(\trf A\trf)
\quad
\mbox{and}\quad
\]

\vspace{-36pt}\vspace{-1pt}
\[
\quad
k\off =\off
k_{\qff \varepsilon}\dff(\trf A\trf)
\off =\off
A
\pff -\pff 
Q_{\trf \varepsilon}\trf(\trf A\trf)
\pff.
\]

\vspace{-12pt}\vspace{3pt}
Clearly,\pss $Q\off =\off f\trf(\trf A\trf)$\sss for some continuous function\sss
$f\dff \colon\dff \rrr\qff \ttoo\qff \rrr$\sss equal\sss to $1$ in a neighborhood of\sss $1$\sss
and\sss to\sss $-\qff 1$ in a neighborhood of\sss $-\qff 1$\nnsp.\oss
Let\sss $q\fff,\qff a$\sss be\sss the images of\sss $Q\fff,\qff A$\sss respectively\sss
in\sss the\dss Calkin\dss algebra of\sss $H$\nnsp.\oss
Then\sss $q\off =\off f\trf(\trf a\trf)$\nnsp.\oss
At\sss the same\sss time\sss the spectrum of\sss $a$\sss is\dss equal\sss to\sss
the essential\sss spectrum of\sss $A$\sss and\sss hence consists of\dss two points\sss
$-\qff 1\fff,\qff 1$\nnsp.\oss
It\sss follows\sss that\sss $a\off =\off f\trf(\trf a\trf)$\sss
and\sss hence\sss $a\off =\off q$\nnsp.\oss
In\sss turn,\oss this implies\sss that\sss $k\off =\off A\qff -\qff Q$\sss is\dss compact.\oss

\myuppar{Compactly-polarized\sss families.}
Let\sss 
$A_{\dff x}\dff \colon\dff
H_{\dff x}\qff \ttoo\qff H_{\dff x}\dff,\pff 
x\qff \in\qff X$\sss
be a\sss family\sss of\dss self-adjoint\sss operators.\oss
We will\sss say\sss that\sss such a family\dss
is\dss a\qss \emph{compactly-polarized\dss family}\oss
if\dss it\dss is\trs fully\trs Fredholm\dss and\sss
all\sss operators\sss $A_{\dff x}$\sss are compactly-polarized.\oss\vspace{-0.27pt}

\mypar{Lemma.}{comp-p-trivializations}
\emph{Suppose\sss that\sss
$A_{\dff x}\dff,\qff x\qff \in\qff X$\sss
is\dss a compactly-polarized\sss family.\oss
Then every\sss strictly adapted\qss ({\fff}local{\fff})\qss trivialization\sss
is\dss fully\sss adapted.\oss}\vspace{-0.27pt}

\proof
Let\sss $z\qff \in\qff X$\nnsp,\qss $\varepsilon\qff \in\qff (\dff 0\fff,\qff 1\dff)$\nnsp,\oss
and\sss $U_{\fff z}$\sss be a neighborhood of\sss $z$\sss
such\sss that\sss $(\trf A_{\dff y}\fff,\qff \varepsilon\trf)$\sss
is\dss an enhanced operator for every\sss 
$y\qff \in\qff U_{\fff z}$\nsp.\oss
Suppose\sss that\sss we are given a strictly adapted\dss 
local\dss trivialization over\sss $U_{\fff z}$\nsp.\oss
Such a\sss local\dss trivialization\sss turns\sss
$P_{\qff [\dff \varepsilon\fff,\qff \infty\dff)}\trf(\trf A_{\dff y}\trf)\fff,\qff
y\qff \in\qff U_{\fff z}$\sss
into a norm continuous family,\oss
and\sss hence also\sss turns\sss
$Q_{\qff \varepsilon}\dff(\trf A_{\dff y}\trf)\fff,\qff
y\qff \in\qff U_{\fff z}$\sss
into a norm continuous family.\oss
It\dss is\dss fully adapted\dss if\trs and\dss only\trs if\trs
it\dss turns\sss
$A_{\dff y}\fff,\qff
y\qff \in\qff U_{\fff z}$\sss
into a norm continuous family.\oss
Since\dss
$A_{\dff y}
\off =\off 
Q_{\qff \varepsilon}\dff(\trf A_{\dff y}\trf)
\qff +\qff 
k_{\qff \varepsilon}\dff(\trf A_{\dff y}\trf)$\nnsp,\oss
for strictly adapted\dss trivializations\sss
the\sss latter condition\dss is\dss equivalent\dss to\sss turning\sss
$k_{\qff \varepsilon}\dff(\trf A_{\dff y}\trf)\fff,\qff
y\qff \in\qff U_{\fff z}$\sss
into a norm continuous family.\oss\vspace{-0.27pt}

Let\sss $K\trf(\trf H\trf)$\sss be\sss the space of\dss
compact\sss operators\sss $H\qff \ttoo\qff H$\sss with\sss the norm\sss topology.\oss
By a\sss theorem of\trs Atiyah\dss and\dss Segal\qss \cite{ase}\qss 
the action of\dss the group\sss
$\mathcal{U}\dff(\trf H\trf)$\sss 
on\sss $K\trf(\trf H\trf)$\sss by conjugations\dss is\dss continuous.\oss
It\dss follows\sss that\dss if\dss the family\sss
$k_{\trf \varepsilon}\trf(\trf A_{\dff y}\trf)\dff,
\off y\qff \in\qff U_{\fff z}$\sss
is\dss norm-continuous\sss in one\sss trivialization,\oss
then\sss it\dss is\dss norm-continuous\sss in every\sss trivialization.\oss\vspace{-0.27pt}

Since\sss $A_{\dff x}\dff,\qff x\qff \in\qff X$\sss is\dss fully\trs Fredholm,\oss
the first\sss paragraph of\dss the proof\dss implies\sss that\sss for some\sss $U_{\fff z}$\sss
some strictly adapted\sss trivialization over\sss $U_{\fff z}$\sss turns\sss 
$k_{\trf \varepsilon}\trf(\trf A_{\dff y}\trf)\dff,
\off y\qff \in\qff U_{\fff z}$\sss
into a norm continuous family.\oss
Then by\sss the previous paragraph\sss every\dss trivialization
over\sss $U_{\fff z}$\sss turns\sss 
$k_{\trf \varepsilon}\trf(\trf A_{\dff y}\trf)\dff,
\off y\qff \in\qff U_{\fff z}$\sss
into a norm continuous family,\oss
and\sss hence every strictly\sss adapted\dss trivialization\sss turns\sss the family\sss
$A_{\dff y}\fff,\qff
y\qff \in\qff U_{\fff z}$\sss
into a norm continuous family.\oss
The\sss lemma\sss follows.\oss  \eproof\vspace{-0.27pt}

\mypar{Theorem.}{index-wss-comp-p}
\emph{Suppose\sss that\sss the family\dss 
$A_{\dff x}\dff,\qff x\qff \in\qff X$\sss is\dss compactly-polarized.\oss
Then a weak spectral\sss section\sss for\sss 
$A_{\dff x}\dff,\qff x\qff \in\qff X$\sss exists\dss if\trs and\dss only\trs if\trs
the analytic\sss index of\qss
$A_{\dff x}\dff,\qff x\qff \in\qff X$\sss vanishes.\oss}\vspace{-0.27pt}

\proof
As in\sss the proof\dss of\trs Theorem\qss \ref{index-wss}\qss
we will\sss assume\sss that\sss $X$\sss is\dss compact.\oss
By\trs Theorem\qss \ref{adapted}\qss there exists a strictly adapted\sss 
trivialization of\sss $\mathbb{H}$\nnsp.\oss
By\dss Lemma\qss \ref{comp-p-trivializations}\qss such\sss trivialization\sss
is\dss fully\sss adapted.\oss
Hence\sss it\dss is\dss sufficient\sss to consider\sss
norm-continuous\sss families\sss of\dss operators\sss 
in\sss a fixed\dss Hilbert\sss space\sss $H$\nnsp.\oss
In\sss this case\trs Theorem\qss \ref{mp}\qss implies\sss
that\sss a weak spectral\sss section exists\sss if\trs and\dss only\trs if\trs
the family\dss is\dss homotopic\sss to a family of\dss invertible operators.\oss
But\sss the\sss latter condition\dss is\dss equivalent\sss to\sss
the vanishing of\dss the analytic index\qss (cf.\qss the proof\dss of\qss
Theorem\qss \ref{index-wss}).\oss  \eproof\vspace{-0.27pt}

\mypar{Theorem}{index-wss-comp-p-full}
\emph{Let\dss 
$A_{\dff x}\dff,\qff x\qff \in\qff X$\sss is\dss 
be a discrete-spectrum\sss and\dss fully\trs Fredholm\dss family.\oss
Then a weak spectral\sss section\sss for\sss  
$A_{\dff x}\dff,\qff x\qff \in\qff X$\sss exists\dss if\trs and\dss only\trs if\trs
its\sss analytic\sss index vanishes.\oss}\vspace{-0.27pt}

\proof
Since\sss $A_{\dff x}\dff,\qff x\qff \in\qff X$\sss is\dss fully\trs Fredholm,\oss
the family\sss 
$C_{\dff x}
\off =\off
\gamma\trf(\trf A_{\dff x}\trf)\dff,\qff x\qff \in\qff X$\sss 
is\dss also fully\trs Fredholm.\oss
Clearly,\oss every operator\sss $C_{\dff x}$\sss is\dss compactly-polarized.\oss
It\sss follows\sss that\sss the family\sss 
$C_{\dff x}\dff,\qff x\qff \in\qff X$\sss is\dss compactly-polarized.\oss
By\trs Theorem\qss \ref{index-wss-comp-p}\qss the family\sss
$C_{\dff x}\dff,\qff x\qff \in\qff X$\sss
admits a weak\sss spectral\sss section\dss if\trs and\dss only\trs if\trs
its analytic index,\oss which\dss is\dss equal\sss to\sss the analytic\sss index of\sss
$A_{\dff x}\dff,\qff x\qff \in\qff X$\nnsp,\oss vanishes.\oss
At\sss the same\sss time every weak\sss spectral\sss section of\dss
$A_{\dff x}\dff,\qff x\qff \in\qff X$\sss
is\dss a weak\sss spectral\sss section of\dss
$C_{\dff x}\dff,\qff x\qff \in\qff X$\nnsp,\oss
and\sss the converse\dss is\dss also\sss true.\oss
The\sss theorem\sss follows.\oss  \eproof

\mypar{Corollary.}{index-ss-alt}
\emph{Suppose\sss that\dss 
$A_{\dff x}\dff,\qff x\qff \in\qff X$\sss is\dss 
a discrete-spectrum\sss and\sss fully\dss Fredholm\dss family.\oss
A spectral\sss section\sss for\sss 
$A_{\dff x}\dff,\qff x\qff \in\qff X$\sss 
exists\dss if\trs and\dss only\trs if\trs the analytic\sss index of\qss
$A_{\dff x}\dff,\qff x\qff \in\qff X$\sss vanishes.}

\proof
It\dss is\dss sufficient\sss to combine\trs
Theorems\qss \ref{index-wss-comp-p-full}\qss 
and\qss \ref{two-types-sections}.\oss  \eproof

\myuppar{Remark.}
While\dss Corollaries\qss \ref{index-ss-alt}\qss and\qss \ref{index-ss}\qss are identical,\oss
the proof\dss of\qss Theorem\qss \ref{index-wss-comp-p-full},\oss 
in contrast\sss with\sss the proof\dss of\trs Theorem\qss \ref{index-wss},\oss
works only\sss for dis\-crete-spectrum\sss families,\oss
because\sss it\dss depends on\sss the properties of\dss compactly-polarized\dss families.\oss

\myuppar{Remark.}
As\sss the reader certainly\sss noticed,\oss we did\sss not\sss really\sss
work with unbounded operators.\oss
The definition of\dss spectral\sss sections\sss looks nicer\sss for unbounded operators,\oss
but\sss can\sss be easily\sss reformulated\sss in\sss terms of\dss their bounded\sss transforms.\oss
Strictly speaking,\oss Corollaries\qss \ref{index-ss-alt}\qss and\qss \ref{index-ss}\qss are
concerned\sss with\sss the image of\dss the bounded\sss transform,\oss
i.e.\qss with\sss the families of\dss self-adjoint\qss 
\emph{strictly\sss contracting}\pss operators.\oss
But\sss in applications such families usually arise from\sss
closed and densely defined\sss unbounded operators.\oss

\myuppar{Remark.}
Recent\sss results of\dss Prokhorova\qss \cite{p3}\qss allow\sss to strengthen\sss
the results of\dss this section.\oss
Namely,\oss in\dss
Lemma\qss \ref{comp-p-trivializations}\qss and\trs Theorem\qss \ref{index-wss-comp-p}\qss
it\dss is\dss sufficient\sss to assume\sss that\sss
$A_{\dff x}\dff,\qff x\qff \in\qff X$\sss is\dss a family of\dss
compactly-polarized operators which\dss is\dss strictly\dss Fredholm\dss and such\sss that\sss
$A_{\dff x}\qff -\qff \lambda\dff,\qff x\qff \in\qff X$\sss
is\dss a\dss Fredholm\dss family\sss for every\sss 
$\lambda\qff \in\qff (\dff -\qff 1 \fff,\qff 1\dff)$\nnsp.\oss
In fact,\oss such families are automatically\sss fully\dss Fredholm.\oss
See\qss \cite{p3},\oss Theorem\qss 5.\oss
In\trs Theorem\qss \ref{index-wss-comp-p-full} 
and\dss Corollary\qss \ref{index-ss-alt}\qss it\dss is\dss sufficient\sss to assume\sss that\sss
$A_{\dff x}\dff,\qff x\qff \in\qff X$\sss is\dss 
a discrete-spectrum\sss and\sss strictly\dss Fredholm\dss family.\oss
As we pointed out\sss at\sss the end of\trs Section\qss \ref{first-proof},\oss
the same\dss is\dss true for\trs Theorem\qss \ref{index-wss}\qss and\dss Corollary\qss \ref{index-ss}.\oss
It\sss seems\sss that\sss strictly\dss Fredholm\dss families provide
a proper context\sss for dealing with spectral\sss sections.\oss

\vspace*{4ex}

\begin{flushright}

November\qss 29,\oss 2021.\oss
Revised\qss January\qss 25,\oss 2022
 
https\halfff:/\!/\!nikolaivivanov.com

E-mail\halfff:\oss nikolai.v.ivanov{\fff}@{\dff}icloud.com

Department\sss of\trs Mathematics,\oss Michigan\sss State\sss University

\end{flushright}


\begin{thebibliography}{MPX}
\vspace{\medskipamount}



\bibitem[A]{a}
M.F.\dss Atiyah,\oss
\emph{$K$\dnsp-theory},\oss
W.A.\dss Benjamin,\oss Inc.,\oss 1967,\oss  iv,\qss 166\qss pp.

\bibitem[ASe]{ase}
M.F.\dss Atiyah,\oss G.\dss Segal,\oss
Twisted $K$\dnsp-theory,\oss
\emph{Ukrainian\dss Mathematical\dss Bulletin},\oss
V.\qss 1,\oss No.\qss 3\qss (2004),\oss 291--334.\oss

\bibitem[{\nsp$AS_{\fff 4}$\dnsp}]{as4}
M.F.\dss Atiyah,\oss I.M.\dss Singer,\oss
The index of\dss elliptic\sss operators\fff:\oss IV,\oss
\emph{Annals of\trs Mathematics},\oss V.\qss 93,\oss No.\qss 1\qss (1971),\oss 
119{\dff}--138.\oss

\bibitem[{\fff\halfff}AS{\halfff}]{as}
M.F.\dss Atiyah,\oss I.M.\dss Singer,\oss
Index\sss theory for skew-adjoint\dss Fredholm\dss operators,\oss
\emph{Publications math\'{e}matiques\qss IHES},\oss 
T.\qss 37\qss (1969),\oss 5{\fff}--26.\oss

\bibitem[DD]{dd}
J.\dss Dixmier,\oss A.\dss Douady,\oss
Champs continus d'espaces hilbertiens et de $C^{*}$\nsp\dnsp-alg\'{e}bres,\oss
\emph{Bulletin de\sss la\sss S.M.F.},\oss T.\qss 91\qss (1963),\oss 227--284.\oss

\bibitem[{\nsp$I_{\dff 1}$\dnsp}]{i1}
N.V.\dss Ivanov,\oss
Topological\sss categories\sss related\sss to\qss Fredholm\qss operators\fff:\oss
I.\oss Classifying\sss spaces,\oss
2021,\oss 107\qss pp.\oss
arXiv:2111.14313.

\bibitem[{\nsp$I_{\dff 2}$\dnsp}]{i2}
N.V.\dss Ivanov,\oss
Topological\sss categories\sss related\sss to\qss Fredholm\qss operators\fff:\oss
II.\oss The analytic\sss index,\oss
2021,\oss 49\qss pp.\oss
arXiv:2111.15081.

\bibitem[K]{ku}
N.\dss Kuiper,\oss
The homotopy\sss type of\dss the unitary\sss group of\trs Hilbert\sss space,\oss
\emph{Topology},\oss
V.\qss 3,\oss No.\qss 1\qss (1964),\oss 19{\dff}--39.\oss

\bibitem[M]{m}
R.\dss Melrose,\oss
Email\sss to\dss M.\dss Prokhorova,\oss August\qss 11,\oss 2021.\oss

\bibitem[{\nsp$MP_{\fff 1}$\dnsp}]{mp1}
R.\dss Melrose,\oss P.\dss Piazza,\oss
Families of\trs Dirac operators,\oss boundaries and\sss the $b$\dnsp-calculus,\oss
\emph{J.\dss of\qss Diff.\dss Geometry},\oss
V.\qss 45\qss (1997),\oss 99{\dff}--180.\oss

\bibitem[{\nsp$MP_{\fff 2}$\dnsp}]{mp2}
R.\dss Melrose,\oss P.\dss Piazza,\oss
An\sss index\sss theorem\sss for\sss families of\trs Dirac\dss operators
on odd-dimensional\sss manifolds with\sss boundary,\oss
\emph{J.\dss of\qss Diff.\dss Geometry},\oss
V.\qss 45\qss (1997),\oss 287--334.

\bibitem[{\nsp$P_{\fff 1}$\dnsp}]{p1}
M.\dss Prokhorova,\oss
Private communication,\oss December\qss 16,\oss 2020.\oss

\bibitem[{\nsp$P_{\fff 2}$\dnsp}]{p2}
M.\dss Prokhorova,\oss
Spectral\sss sections,\oss 2020 -- 2021,\oss 37\qss pp.\oss
arXiv:2008.04672v5,\oss
\emph{Israel\dss Journal\dss of\dss Mathematics},\oss
to appear.\oss

\bibitem[{\nsp$P_{\fff 3}$\dnsp}]{p3}
M.\dss Prokhorova,\oss
The continuity properties of discrete-spectrum families of Fredholm operators,\oss 2022,\oss
4\qss pp.\qss arXiv:2201.09869.

\bibitem[Se]{see}
R.\dss Seeley,\oss
Complex\sss powers of\dss an elliptic operator,\oss
\emph{Proceeding of\dss Symposia in\sss Pure\sss Mathematics},\oss
V.\qss 10,\oss American\sss Mathematical\sss Society,\oss 1967,\oss
pp.\qss 288--307.\oss

\bibitem[S]{sfc}
G.\dss Segal,\oss
Fredholm\dss complexes,\oss
\emph{Quarterly\trs J.\dss of\trs Mathematics},\oss
V.\qss 21\qff (1970),\oss 385--402.



\end{thebibliography}
\end{document}